\documentclass[11pt]{article}
\usepackage{authblk}
\usepackage{amsfonts,amsmath,amssymb,amsthm}
\usepackage{latexsym}
\usepackage{mathrsfs}
\usepackage{tikz}
\usepackage{color}
\usepackage{multirow}

\newcommand{\uc}{U_q(\mathfrak{su}(2))}
\newcommand{\unc}{U_q(\mathfrak{su}(1,1))}

\newcommand{\luc}{\mathfrak{su}(2)}
\newcommand{\lunc}{\mathfrak{su}(1,1)}
\newcommand{\osc}{\mathfrak{osc}}

\newcommand{\fa}{\mathfrak{a}}
\newcommand{\fc}{\mathfrak{c}}
\newcommand{\fe}{\mathfrak{e}}
\newcommand{\ff}{\mathfrak{f}}
\newcommand{\fh}{\mathfrak{h}}
\newcommand{\fn}{\mathfrak{n}}

\newcommand{\CC}{\mathbb{C}}

\newcommand{\cX}{\mathcal{X}}
\newcommand{\cY}{\mathcal{Y}}

\begin{document}

\title{New realizations of algebras of the Askey--Wilson type\\
in terms of Lie and quantum algebras.}

\author[$\dagger$]{Nicolas Cramp\'e}

\author[$\star$]{Dounia Shaaban Kabakibo}

\author[$\star$]{Luc Vinet}

\affil[$\dagger$]{Institut Denis-Poisson CNRS/UMR 7013 - Universit\'e de Tours - Universit\'e d'Orl\'eans,
Parc de Grandmont, 37200 Tours, France.}

\affil[$\star$]{Centre de recherches math\'ematiques, Universit\'e de Montr\'eal, P.O. Box 6128, 
Centre-ville Station, Montr\'eal (Qu\'ebec), H3C 3J7, Canada.}

\maketitle

\begin{abstract}

The Askey--Wilson algebra is realized in terms of the elements of the quantum algebras $\uc$ or $\unc$.
A new realization of the Racah algebra in terms of the Lie algebras $\luc$ or $\lunc$ is given also.
Details for different specializations are provided.
The advantage of these new realizations is that one generator of the Askey--Wilson (or Racah) algebra becomes diagonal in the usual representation of the quantum algebras 
whereas the second one is tridiagonal.
This allows to recover easily the recurrence relations of the associated orthogonal polynomials of the Askey scheme.
These realizations involve rational functions of the Cartan generator of the quantum algebras, they are linear with respect to the other generators and depend on the Casimir element of the quantum algebras.

\end{abstract}

%

%


\section{Introduction}

The Askey--Wilson algebra has been introduced in \cite{Z91} as the algebra realized by the recurrence operator $\cY$
and the $q$-difference operator $\cX$ intervening in the bispectral problem associated to the Askey--Wilson polynomials (see for example \cite{KLS} 
for an introduction to these orthogonal polynomials). This explains the name of the algebra which also appears
in different other contexts. It has been used in the framework of integrable systems as a quotient of the reflection algebra \cite{Bas}.
It is connected with the symmetry of the $6j$ or Racah coefficients of $\unc$ \cite{GZ,GZ2} as the centralizer of the 
diagonal action of this quantum algebra in its three-fold tensor product (see also \cite{CPV} for recent advances in this subject).
Let us mention that similar connections exist for other pairs of algebras: replace the couple Askey--Wilson algebra and $\unc$ by $q$-Racah algebra and $\uc$ or by Racah algebra and $\luc$.
There exists another connection between the Askey--Wilson algebra and $\unc$. Indeed, the two operators $\cX$ and $\cY$ can be realized as
combinations of the raising and lowering generators of $\unc$ \cite{GZ93}. 
This realization allows to use the irreducible representations of $\unc$ to obtain representations of the Askey--Wilson algebra.
In this case, the operators $\cX$ and $\cY$ become tri-diagonal matrices and 
the overlap functions between eigenstates of these matrices are expressed via orthogonal polynomials of the Askey scheme \cite{KV,BMVZ,AAK}.

In the realizations obtained up to now of the Askey--Wilson algebra using twisted primitive elements \cite{K,GZ93,BMVZ} or using equitable representations \cite{ITW,Ter}, 
both operators $\cX$ and $\cY$ are tri-diagonal or bi-diagonal in the usual representations of $\uc$ or $\unc$. Then, one of these operators (or their $q$-commutator) is diagonalized 
in some basis via the help of orthogonal polynomial and the other ones remains tridiagonal. At this stage, one can recognize, the recurrence operators and 
the $q$-difference operators defining orthogonal polynomials.
In this paper, we give a new realization of the Askey--Wilson algebra in terms of elements in $\unc$ such that $\cX$ (resp. $\cY$) are immediately diagonal (resp. tridiagonal)
in the standard representation of $\unc$.
To do so, we look for a realization of $\cX$ and $\cY$ of the following form
\begin{eqnarray}
 \cX= f_1(K) \quad \text{and} \qquad \cY=E+f_2(K, C) + f_3(K, C) F\ , \label{eq:real}
\end{eqnarray}
where $E$, $K$ and $F$ are the generators of $\unc$ and $C$ is its quadratic Casimir element (see Subsection \ref{def:ql} for the exact definition). 
The functions $f_1$, $f_2$ and $f_3$ are to be determined so that $\cX$ and $\cY$ satisfy the defining relations of the Askey--Wilson algebra. In the usual representation of $\unc$, $K$ is represented by 
a diagonal matrix therefore $\cX$ is also diagonal. Let us emphasize that to succeed in finding such a realization, it is necessary that $f_2$ and $f_3$ depend on the Casimir element and
that they be rational functions. To be precise, because we need rational functions in terms of the Cartan generator, we take the inverse of certain polynomials in the Cartan generator 
(see for example \eqref{eq:AWr}). 
To do so, we must use a localization of $\unc$ where these elements are invertible. In the rest of this paper, by abuse of language, we shall use the term quantum algebra and the notation $\unc$ for this type of localization. 

In Section \ref{sec:nd}, we focus on the non-deformed case. We find a new realization of the Racah algebra in terms of the generators of $\luc$ or $\lunc$ (see Subsection \ref{sec:racah}). We then prove that using 
the finite irreducible representations of $\luc$, we recover directly the recurrence operator associated to the Racah polynomials whereas using the infinite irreducible representations of $\lunc$, 
we obtain the one of the Wilson polynomials (see Subsection \ref{sec:racah}). One thus arrives at the polynomials sitting at the top of the Askey scheme.
In Subsection \ref{sec:spec}, we provide the specializations or contractions of the realization proposed so as to obtain the recurrence operators associated to all the polynomials of the Askey scheme.
For the polynomials at the bottom of the tableau, we must replace the Lie algebra $\lunc$ by the algebra of the oscillator $\osc$ (which can be seen as a contraction of $\lunc$).

In Section \ref{sec:defo}, we deal with the $q$-deformed case and provide the explicit realization of the Askey--Wilson (or $q$-Racah) algebra of the form \eqref{eq:real}.
In Subsection \ref{sec:AW}, we prove that this realization allows to obtain the recurrence operator of $q-$Racah polynomials by using finite representations of $\uc$. 
As in the non-deformed case, we give different specializations and contractions of this realization that corresponds to various families of $q$-orthogonal polynomials.

\section{Racah algebra, its specializations and polynomials of the Askey scheme \label{sec:nd}}

The Askey--Wilson algebra $AW(3)$ is generated by two elements $X$ and $Y$ obeying the following defining equations \cite{Z91}
\begin{eqnarray}
{}&&  X^2Y-\beta XYX+YX^2 =\gamma(XY+YX)+\gamma^*X^2+ \omega X + \rho Y +\eta , \label{eq:AW1} \\
{}&& Y^2X-\beta YXY+XY^2=\gamma^*(XY+YX)+\gamma Y^2+ \omega Y + \rho^* X +\eta^*, \label{eq:AW2}
\end{eqnarray}
where $\beta$, $\gamma$, $\gamma^*$, $\omega$, $\rho$, $\eta$, $\rho^*$ and $\eta^*$ are parameters. 

In this section, we focus on the case $\beta=2$. The algebra is usually called Racah algebra and has been introduced 
in \cite{GZ2} (see \cite{GVZ} for a review) to study the Racah $W$-coefficients or $6j$-symbols intervening in the coupling of three angular momenta.
We give a new realization of the Racah algebra and some specializations in terms of $\luc$, $\lunc$ or $\osc$.
As explained in the introduction, we want that one of the generators of the Racah algebra be diagonal in the standard representation of these Lie algebras.
Therefore, we look for a realization which takes the following form 
\begin{eqnarray}
 X= f_1(\fh) \quad \text{and} \qquad Y=\fe+f_2(\fh,\fc) + f_3(\fh,\fc) \ff \ ,\label{eq:real-lie}
\end{eqnarray}
where $\fe$, $\ff$ and $\fh$ are the usual raising, lowering and Cartan generators of $\lunc$ or $\luc$. Since $\fh$ is diagonal in the representation of $\luc$ or $\lunc$, we obtain that $X$ is diagonal while 
$Y$ becomes the recurrence operator defining orthogonal polynomials of the $q=1$ Askey scheme.
We recover the whole list of polynomials given in \cite{KLS}. In Table \ref{tw}, we summarize the different specializations 
of the Racah algebra we obtained with the name of the associated algebras (the parameters refer to those of equations \eqref{eq:AW1}-\eqref{eq:AW2}).
\begin{table}[htb]
\begin{center}
 \begin{tabular}{|c|c|c|c|c|c|c|c|c|c|}
  \hline
  $\beta$ & $\gamma$ &$\gamma^*$ & $\omega$ & $\rho$ & $\eta$ & $\rho^*$ & $\eta^*$ & Algebras\\
  \hline
\multirow{2}{*}2& \multirow{2}{*}{2}& \multirow{2}{*}{2} & \multirow{2}{*}{$2(d+\fc)$} & \multirow{2}{*}{$a^2-1$} & \multirow{2}{*}{$2d\fc$} 
        &  $4\fc-1 $  & $(b+c+1)\fc$ & Racah  \\
 & & & & & &$+(c-b)^2$ &$\times (2a-b-c-1)$ &  \\
  
  \hline
 2 & 2 & 0 & $\alpha-2\beta-1 $ &$\alpha^2-1$&$( \alpha-2\beta-1) \fc$&$1$&$-\fc$&Hahn  \\
 
   \hline
 2 & 2 & 0 & 0 &$\alpha^2-1$&$\displaystyle\frac{1}{2}(1-\alpha^2)+2 \fc$&$ 0 $&$\displaystyle -\frac{1}{2}$&Jacobi  \\

  \hline
 2 & 0 & 2 & $-1-\mu-\nu$ &$1$&$0$&$ (\mu+\nu)^2-1+4\mathfrak{c}  $&$-2(1+\mu+\nu)\fc$&Dual Hahn  \\

  \hline
 2 &0&0&$b$&$1$&$0$&$b^2\pm 4$&$0$&Lie type \\
  \hline
 2& {0}& {$-1$} & {$-b$} & {1} & {0} &   {$b^2$} &{$-2$} & Oscillator   \\
  \hline
\end{tabular}
\caption{Different specializations of the Racah algebra and associated polynomials 
($d$ is given in \eqref{eq:delta} and $\fc$ is the Casimir element of $\luc$ or $\lunc$, the upper (resp. lower) sign is for $\luc$ (resp. $\lunc$) ).\label{tw}}
\end{center}
\end{table}

\subsection{Definitions of the Lie algebras $\luc$, $\lunc$ and of the oscillator algebra $\osc$ \label{sec:def}}

In this section, we recall well-known results about the Lie algebras $\luc$, $\lunc$ and $\osc$ as well as their relevant representation 
to introduce the notations used in this paper.

The Lie algebras $\luc$ and $\lunc$ are generated by $\fe$, $\ff$ and $\fh$  subject to
\begin{eqnarray}
 [\fh,\fe]=\fe \ , \qquad [\fh,\ff]=-\ff  \ ,\qquad [\fe,\ff]=\pm 2\fh\ . \label{eq:ldefu}
\end{eqnarray}
The upper (resp. lower) sign in the previous relation corresponds to $\luc$ (resp. $\lunc$).
The Casimir element is given by 
\begin{equation}
 \mathfrak{c}=   \pm \fe\ff + \fh(\fh-1)\ .
\end{equation}
The oscillator algebra $\osc$ is generated by $\fa$ and $\fa^{\dagger}$ subject to
\begin{equation}
 [\fa,\fa^\dagger]=1\ .
\end{equation}
We define as usual the number operator by $\fn=\fa^\dagger\fa$ such that $[\fn,\fa]=-\fa$ and $[\fn,\fa^\dagger]=\fa^\dagger$.

The finite irreducible representations for $\luc$ are characterized by an integer or half-integer $j$ and are given explicitly by
\begin{eqnarray}
 \fe |n\rangle  &=& |n-1\rangle \quad\text{for } 1\leq n\leq 2j \quad\text{and} \qquad  \fe |0\rangle=0,\\
 \ff |n\rangle  &=& (2j-n)(n+1)\ |n+1\rangle \quad\text{for } 0\leq n\leq 2j ,\\
 \fh |n\rangle  &=& (j-n)\  |n\rangle \quad\text{for } 0\leq n\leq 2j,
\end{eqnarray}
where $\{ \ |n\rangle \ |\  0\leq n \leq 2j \ \}$ is a basis of $\CC^{2j+1}$.
In this representation, the Casimir element $\fc$ is $j(j+1)$ times the identity matrix.

The infinite irreducible representations for $\lunc$ of the positive discrete series are characterized by a positive real $\ell$ and are given explicitly by
\begin{eqnarray}
 \fe |n\rangle  &=& |n-1\rangle \quad\text{for }  n\geq 1 \quad\text{and} \qquad  \fe |0\rangle=0,\\
 \ff |n\rangle  &=& (2\ell+n)(n+1)\ |n+1\rangle \quad\text{for }  n\geq 0, \\
 \fh |n\rangle  &=& -(\ell+n)\  |n\rangle \quad\text{for }  n \geq 0.
\end{eqnarray}
In this representation, the Casimir element $\fc$ is $\ell(\ell-1)$ times the identity matrix.

The infinite irreducible representations for $\osc$ are given by
\begin{eqnarray}
 \fa |n\rangle  &=& |n-1\rangle \quad\text{for }  n\geq 1 \quad\text{and} \qquad  \fa |0\rangle=0\ ,\\
 \fa^\dagger |n\rangle  &=& (n+1)\ |n+1\rangle \quad\text{for }  n\geq 0. \\
\end{eqnarray}

\subsection{Racah algebra and Wilson/Racah Polynomials \label{sec:racah}}

In this subsection, we give an explicit form of the operators $X$ and $Y$ in terms of the generators of $\uc$ and $\unc$ so that they satisfy the Racah algebra relations.
Indeed, if we define these operators by
\begin{eqnarray}
X&=& \fh(\fh-a )\ , \\
Y&=&\fe-2 \frac{(\fh^2-a\fh+\fc)(\fh^2-a\fh+d) }{(2\fh-a+1)(2\fh-a-1)} \nonumber\\
&& \mp  \frac{(\fh^2+(1-2a)\fh-\fc+(a-1)a )(\fh-b)(\fh-a+b+1)(\fh-c)(\fh-a+c+1)}
{(2\fh-a)(2\fh-a+1)^2(2\fh-a+2)}\ \ff \ , \label{eq:Yracah}
\end{eqnarray}
with the upper (resp. lower) sign in $Y$ for $\luc$ (resp. $\lunc$) and
\begin{equation}
 d=\frac{1}{2}(a+1)(b+c+1)-(b+1)(c+1)\ ,\label{eq:delta}
\end{equation}
we can show by direct computations that they satisfy the relations
\begin{eqnarray}
&&  [X,[X,Y]] =2(XY+YX)+2 X^2+ 2(d+\fc) X + (a^2-1) Y +2 d \fc , \\
&&  [Y,[Y,X]]= 2(XY+YX)+2 Y^2+ 2(d+\fc) Y +(4\fc+(c-b)^2-1 ) X \nonumber\\
&& \hspace{5cm}+(b+c+1)(2a-b-c-1)\fc  ,
\end{eqnarray}
which are recognized as those of the Racah algebra.
Let us remark that we use the following exchange relations to perform the computations 
\begin{equation}
 \fe   \frac{1}{2\fh-a} =  \frac{1}{2\fh-a-2}  \fe \quad \text{and} \qquad  \ff   \frac{1}{2\fh-a} =  \frac{1}{2\fh-a+2}  \ff  \ .
\end{equation}
We have introduced three free parameters $a$, $b$ and $c$ in this realization of the Racah algebra. 

The operator $X$ (resp. $Y$) given by \eqref{eq:real-lie} in the finite irreducible representation of $\luc$ becomes a diagonal (resp. tridiagonal) matrix. 
Therefore, if we set $\displaystyle |p\rangle=\sum_{n=0}^{2j} p_n |n\rangle$, the eigenvalue problem
\begin{equation}
 Y|p\rangle =\lambda |p\rangle \label{eq:vp}
\end{equation}
gives the following recurrence relation for the components of the vector: 
\begin{equation}
 p_{n+1}+ f_2\big(j-n,j(j+1)\big) p_n +n(2j-n+1)\ f_3\big(j-n,j(j+1)\big)  p_{n-1}=\lambda p_n\ , \label{eq:recg1}
\end{equation}
with the convention $p_{2j+1}=p_{-1}=0$. 
Similarly for $\lunc$, the following recurrence relation is obtained: 
\begin{equation}
 p_{n+1}+ f_2\big(-\ell-n,\ell(\ell-1)\big) p_n +n(2\ell+n-1)\ f_3\big(-\ell-n,\ell(\ell-1)\big)  p_{n-1}=\lambda p_n\ . \label{eq:recg2}
\end{equation}

When $f_2$ and $f_3$ are given as per \eqref{eq:Yracah}, these recurrence relations become the ones of the Racah or Wilson polynomials.
More precisely, one finds that the recurrence relations are the ones of 
\begin{itemize}
\item  the Racah polynomial (see Section 1.2 of \cite{KLS}) for the
Lie algebra $\luc$, with $x=0,\dots, 2j$,
\begin{equation}
 p_n=\frac{(c-j+1)_n(a-b-j)_n(-2j)_n}{(n-2j+a)_n}
R_n(\tau(x);c-j,a-c-j-1,-2j-1,c-b) ,
\end{equation}
with $\tau(x)=x(x-2j+c-b)$ and $\lambda=\tau(x)+j(1-b+c)$;
\item the Wilson polynomial  (see Section 1.1 of \cite{KLS}) for the
Lie algebra $\lunc$,
 \begin{equation}
 p_n=\frac{1}{(n+a+2\ell)_n}W_n\left(x^2;a-\frac{b+c}{2},\ell-\frac{b-c}{2},\ell+\frac{b-c}{2},1+\frac{b+c}{2}\right)
,
\end{equation}
with $\lambda=\frac{(b-c)^2}{4}+\ell(\ell-1)-x^2$.
\end{itemize}
We have collected in the first row of the Table \ref{tw} the results
of this subsection.

\subsection{Specializations of the Racah algebra and polynomials of the Askey scheme \label{sec:spec}}

We study now specializations or contractions of the realization of the Racah algebra given by \eqref{eq:Yracah}.
We show that we obtain the whole list of the recurrence relations of the polynomials of the Askey scheme.
The number of free parameters decreases as the algebra is specialized. 

\paragraph{Hahn algebra.} 

Let us introduce a realization of the Racah algebra depending on two free parameters.  
The following two operators
\begin{eqnarray}
X&=& \fh(\fh-\alpha )\ , \\
Y&=&\fe-( \alpha-2\beta-1) \frac{\fh^2-\alpha\fh+\fc }{(2\fh-\alpha+1)(2\fh-\alpha-1)} \nonumber\\&& 
\pm   \frac{(\fh^2+(1-2\alpha)\fh-\fc+(\alpha-1)\alpha )(\fh-\beta)(\fh-\alpha+\beta+1)}{(2\fh-\alpha)(2\fh-\alpha+1)^2(2\fh-\alpha+2)}\ \ff  \ ,
\end{eqnarray}
with the upper sign in $Y$ for $\luc$ and the lower sign for $\lunc$, satisfy the specialization of the Racah algebra relations which define the Hahn algebra.
We give in Table \ref{tw} the values of the parameters of relations \eqref{eq:AW1}-\eqref{eq:AW2} in this case.

The diagonalization of $Y$ are associated to 
\begin{itemize}
 \item the Hahn polynomial (see Section 1.5 of \cite{KLS}) for the Lie algebra $\luc$, with $x=0,\dots, 2j$,
 \begin{equation}
  p_n=(-1)^n\frac{(\alpha+1)_n(-2j)_n}{(n+\alpha-2j)_n} Q_n(x;\beta-j,\alpha-\beta-j-1,2j) , \qquad  \lambda=j-x \ ;
 \end{equation}
 \item the continuous Hahn polynomial (see Section 1.4 of \cite{KLS}) for the Lie algebra $\lunc$ 
\begin{equation}
  p_n=\frac{n!}{i^n(n+2\ell+\alpha)_n}p_n(x;\beta+\ell-d+1,2\ell-d,d-\beta-\ell+\alpha,d) , \qquad  \lambda=-ix+d-\ell.
 \end{equation}
\end{itemize}

Let us remark that the Hahn algebra is isomorphic to the Higgs algebra \cite{Z92,LV,GVZ,FGVVZ} and to a truncation of the reflection equation associated to a Yangian \cite{CRVZ}.

\paragraph{Dual Hahn algebra.} 

The following operators provide a two-parameter realization of the Racah algebra  
\begin{eqnarray}
X&=& \fh ,\\
Y&=&\fe-2\fh^2+(1+\mu+\nu)\fh \mp (\fh-\mu)(\fh-\nu)\ff ,
\end{eqnarray}
with the upper sign in $Y$ for $\luc$ and the lower sign for $\lunc$.
The algebra generated in this case is called the dual Hahn algebra and is isomorphic to the Hahn algebra.
The structure constants are again in Table \ref{tw}.

The eigenvectors diagonalizing $Y$ are associated to the following polynomials:  
\begin{itemize}
 \item  for $\luc$, the dual Hahn polynomial (see Section 1.6 of \cite{KLS})
 \begin{equation}
  p_n=(\mu-j+1)_n (-2j)_n R_n(\tau(x);\mu-j,-\nu-j-1,2j), \quad, \quad \lambda=\tau(x)+j(\nu-\mu-1)\ ,
 \end{equation}
 where $ \tau(x)=x(x-2j+\mu-\nu)$;
 \item  for $\lunc$, the continuous dual Hahn polynomial (see Section 1.3 of \cite{KLS})
 \begin{equation}
  p_n=S_n\left(x^2;\ell+\frac{\mu-\nu}{2},\ell-\frac{\mu-\nu}{2},1+\frac{\mu+\nu}{2}\right), \quad \lambda=-x^2-\ell(1-\ell)-\frac{(\mu-\nu)^2}{4} \ .
 \end{equation}
\end{itemize}

\paragraph{Jacobi algebra.} 

The following operators, depending on one free parameter, 
\begin{eqnarray}
X&=& \fh(\fh-\alpha )\ , \\
Y&=&\fe- \frac{1-\alpha^2+  4\fc }{2(2\fh-\alpha+1)(2\fh-\alpha-1)} \mp \frac{\fh^2+(1-2\alpha)\fh-\fc+(\alpha-1)\alpha }{(2\fh-\alpha)(2\fh-\alpha+1)^2(2\fh-\alpha+2)}\ \ff \ ,
\end{eqnarray}
with the upper sign in $Y$ for $\luc$ and the lower sign for $\lunc$, realize a specialization of the Racah algebra relations which pertain to the Jacobi algebra.
The structure constants of this algebra are listed in Table \ref{tw}.

In this case, the recurrence relation, for $\lunc$ and $\alpha \in ]-1,\infty[$, are those of the Jacobi polynomials  (see Section 1.8 in \cite{KLS}) where $\ell= \frac{\beta+1 }{2}$ with 
$\beta \in  ]-1,\infty[$ and
 \begin{equation}
p_n(x)= \frac{n!}{n+\alpha+\beta +1} P_n ^{(\alpha ,\beta)} (x) \quad \lambda= \frac{x}{2}   \quad \text{for } x \in [-1,1] \ .
\end{equation}

\paragraph{Lie type algebras.}

The following linear combinations of the generators of the Lie algebras $\luc$ or $\lunc$ 
\begin{eqnarray}
X&=& \fh ,\\
Y&=&\fe-b\fh + \ff ,\label{eq:YKRAW}
\end{eqnarray}
lead also to a specialization of the Racah algebra. In this case, the Racah algebra becomes a Lie algebra.
The explicit specialization of the parameters of the Racah algebra is given in Table \ref{tw}.

The eigenvectors diagonalizing $Y$ are associated to the following families of polynomials:
\begin{itemize}
 \item for $\luc$ and $b \in  {\rm \mathbb{R}}$, one obtains the Krawtchouk polynomials (see Section 1.10 of \cite{KLS}) by letting $b=\epsilon \frac{ (1-2p)}{\sqrt{p(1-p)}}$ with $\epsilon  \in \{1,-1\}$, $p \in  ]0,1[$ and
 \begin{equation}
  p_n(x) =\left( \frac{\epsilon}{\sqrt{p(1-p)}} \right)^n (-2j)_np^nK_n(x;p,2j),  \quad \lambda= \frac{\epsilon(x-j)}{\sqrt{p(1-p)}}   \quad \text{for } x=0,1,\dots 2j\ ;
 \end{equation}
\item for $\lunc$ and
\begin{itemize}
\item $b \in ]-\infty,-2[  \cup ]2,\infty[$, one finds the Meixner polynomials (see Section 1.9 of \cite{KLS}) by letting $b=\epsilon \left(\frac{1}{{c}}+{c}\right)$ with $\epsilon\in  \{1,-1\}$, $c\in  ]0,1[$ and
\begin{equation}
  p_n(x) = ( 2\ell)_n\left(\frac{-\epsilon c}{c+1}\right)^nM_n(x;2\ell,c^2), \quad \lambda=\epsilon \left(\frac{1}{{c}}-{c}\right)(x+\ell) \quad \text{for } x=0,1,\dots  \ ;
 \end{equation}
 \item $b \in ]-2,2[ $, one gets the Meixner-Pollaczek polynomials (see Section 1.7 of \cite{KLS}) by letting $b=-2 \epsilon \cos{\phi}$ with $\epsilon  \in \{1,-1\}$, $\phi \in ]0,\pi[$ and 
 \begin{equation}
  p_n(x) =\epsilon^n n! P_n ^{(\ell)}(x;\phi), \quad \lambda=2 \epsilon  \sin{(\phi)} x \quad \text{for }  x \in \mathbb{R} \ ;
   \end{equation}
 \item $b \in   \{-2,2\} $, one arrives at the Laguerre polynomials (see Section 1.11 of \cite{KLS}) by letting $b=2\epsilon $ with $\epsilon  \in \{1,-1\} $  and
  \begin{equation}
    p_n(x) =(-\epsilon)^n n! L_n^{(2\ell-1)}(x) ,   \quad \lambda=\epsilon x \quad \text{for }x \in [0,\infty[\ .
 \end{equation}
\end{itemize}
\end{itemize}

\paragraph{Oscillator Racah algebra.}

Finally, we give here a realization of the Racah algebra in terms of the oscillator algebra $\osc$.
The following linear combinations of the generators of the oscillator algebra $\osc$
\begin{eqnarray}
X&=& \fn ,\\
Y&=&\fa +b\ \fn  + \fa^\dagger. \label{eq:Yosc}
\end{eqnarray}
generate a specialization of the Racah algebra called the oscillator Racah algebra.
The parameters are again given in Table \ref{tw}.

The associated polynomials are in this case:
\begin{itemize}
\item  $b \in \mathbb R_{\ne 0}$, one gets the Charlier polynomials  (see Section 1.12 in \cite{KLS}) by letting $b=\frac{\epsilon}{\sqrt{a}}$ with a $>$ 0 and 
 \begin{equation}
   p_n(x) =(-\epsilon \sqrt{a})^{n}C_n(x;a) ,   \quad \lambda=\epsilon \frac{x-a}{\sqrt{a}} \quad \text{for } x=0,1,\dots \ ;
 \end{equation}
 \item $b=0$, one finds the Hermite polynomials  (see Section 1.13 in \cite{KLS}) by letting 
 \begin{equation}
   p_n(x) =2^{-\frac{n}{2}}H_n(x) ,   \quad \lambda=\sqrt{2}x \quad \text{for }x \in  \mathbb{R} \ .
\end{equation}
\end{itemize}

\section{Deformed case \label{sec:defo}}

We study now the case $\beta\neq \pm2$, that is $q$ not equal to $1$.
After a shift of the generators and by setting $\beta=q^2+q^{-2}$, the defining relations \eqref{eq:AW1}-\eqref{eq:AW2} of $AW(3)$ can be rewritten in the form \cite{Vid, Vid2}
\begin{eqnarray}
{} [\cX,[\cX,\cY]_q]_{q^{-1}} &=& \omega \cX + \rho \cY +\eta, \label{eq:AW3} \\
{}[\cY,[\cY,\cX]_q]_{q^{-1}} &=& \omega \cY + \rho^* \cX +\eta^*,\label{eq:AW4} 
\end{eqnarray}
where $[\cX,\cY]_q=q\cX\cY-q^{-1}\cY\cX$. We suppose that $q$ is not a root of unity.

We shall give different realizations of the Askey--Wilson algebra for $\beta \neq \pm 2$ in terms of $\uc$ or $\unc$.
In the irreducible representation of $\uc$ or $\unc$, $\cY$ will turn into the difference operator that defines the orthogonal polynomials of the $q$-deformed part of the Askey scheme \cite{KLS}.
We can recover in this way all the different polynomials given in \cite{KLS}, we shall however limit ourselves to providing some explicit examples.

\subsection{Quantum Lie algebras $\uc$ and $\unc$ \label{def:ql}}

We here recall the definition of the quantum algebras $\uc$ and $\unc$ and also introduce (some of) their irreducible representations.

The quantum Lie algebras $\uc$ and $\unc$ are generated by $E$, $F$ and $q^{H}$  subject to
\begin{eqnarray}
 q^H E=qE q^H \ , \qquad q q^H  F=Fq^H  \ ,\qquad [E,F]=\pm [2H]_q\ , \label{eq:defu}
\end{eqnarray}
where $[X]_q=\frac{q^X-q^{-X}}{q-q^{-1}}$.
The upper (resp. lower) sign in the previous relation corresponds to $\uc$ (resp. $\unc$).
The Casimir element is given by 
\begin{equation}
 C= (q-q^{-1})^2( \pm EF + [H]_q\ [H-1]_q)+q+q^{-1}\ . \label{eq:Cb}
\end{equation}

The finite irreducible representations for $\uc$ are characterized by an integer or a half-integer $j$ and are given explicitly by
\begin{eqnarray}
 E |n\rangle  &=& |n-1\rangle \quad\text{for } 1\leq n\leq 2j \quad\text{and} \qquad  E |0\rangle=0\ ,\\
 F |n\rangle  &=& [2j-n]_q[n+1]_q\ |n+1\rangle \quad\text{for } 0\leq n\leq 2j \ ,\\
 K |n\rangle  &=& q^{j-n}\  |n\rangle \quad\text{for } 0\leq n\leq 2j\ ,
\end{eqnarray}
where $\{ \ |n\rangle \ |\  0\leq n \leq 2j \ \}$ is a basis of $\CC^{2j+1}$.
In this representation the Casimir element $C$ is $q^{2j+1}+q^{-2j-1}$ times the identity matrix.

One infinite irreducible representations series for $\unc$ is characterized by a positive real $\ell$ and given explicitly by
\begin{eqnarray}
 E |n\rangle  &=& |n-1\rangle \quad\text{for }  n\geq 1 \quad\text{and} \qquad  E |0\rangle=0\ ,\\
 F |n\rangle  &=& [2\ell+n]_q[n+1]_q\ |n+1\rangle \quad\text{for }  n\geq 0\ , \\
 K |n\rangle  &=& q^{-\ell-n}\  |n\rangle \quad\text{for }  n \geq 0 \ .
\end{eqnarray}
In this representation, the Casimir element $C$ is $q^{2\ell-1}+q^{-2\ell+1}$ times the identity matrix.

\subsection{Askey--Wilson algebra and Askey--Wilson/$q$-Racah polynomials \label{sec:AW} }

In the following, we give a three-parameter realization of the operators $\cX$ and $\cY$ in the form \eqref{eq:real} which satisfies the Askey--Wilson algebra relations. 
These operators $\cX$ and $\cY$ are (the upper (resp. lower) sign is for $\uc$ (resp. $\unc$) ):
\begin{eqnarray}
\cX&=& K^2-aK^{-2}\ , \\
\cY&=&E+ K^2 \frac{\left(\frac{(1-a)(a-bc)}{a}+(b+c)  C\right)(K^4-a)+\left((b+c)(a-1)+(a-bc)  {C}\right) (q+q^{-1})K^2   }{(q-q^{-1})(q^{-2}K^4+a)(q^{2}K^4+a)} \nonumber \\
&& \pm q K^2\  \frac{(a^2+aq  {C}K^2+q^2K^4)(bqK^2+a)(qK^2-b)(qK^2-c)(qc K^2 +a)}{a(K^4+a)(q^2K^4+a)^2(q^{4}K^4+a)}  F \ ,  \label{eq:AWr}
\end{eqnarray}
and they satisfy the relations \eqref{eq:AW3}-\eqref{eq:AW4} of $AW(3)$ that read here:
\begin{eqnarray}
{} \frac{1}{q^2-q^{-2}}[\cX,[\cX,\cY]_q]_{q^{-1}} &=& \frac{1}{q+q^{-1}}\left(\frac{(a-1)(a-bc)}{a}-(b+c)  C  \right) \cX + a(q^2-q^{-2}) \cY \nonumber\\
&& +(1-a)(c+b)-(a-bc) {C}, \label{eq:AW3g} \\
{}\frac{1}{q+q^{-1}}[\cY,[\cY,\cX]_q]_{q^{-1}} &=& \frac{q-q^{-1}}{q+q^{-1}}\left(\frac{(a-1)(a-bc)}{a}-(b+c)  C  \right) \cY + \frac{bc}{a}(q+q^{-1}) \cX\nonumber\\
&&+\frac{(a-bc)(b+c)}{a}+   \frac{(a-1)bc}{a} {C}   .\label{eq:AW4g} 
\end{eqnarray}

As explained in Section \ref{sec:racah}, the components of the eigenvectors diagonalizing $\cY$ of the form \eqref{eq:real} satisfy a three-term recurrence relation.
Here, in the deformed case, for $\uc$, it is
\begin{equation}
 p_{n+1}+ f_2\big(q^{j-n},q^{2j+1}+q^{-2j-1}\big) p_n +[2j-n+1]_q[n]_q\ f_3\big(q^{j-n},q^{2j+1}+q^{-2j-1}\big)  p_{n-1}=\lambda p_n\ , \label{eq:recg3}
\end{equation}
and for $\unc$,
\begin{equation}
 p_{n+1}+ f_2\big(q^{-\ell-n},q^{2\ell-1}+q^{-2\ell+1}\big) p_n +[2\ell+n-1]_q[n]_q\ f_3\big(q^{-\ell-n},q^{2\ell-1}+q^{-2\ell+1}\big)  p_{n-1}=\lambda p_n\ . \label{eq:recg4}
\end{equation}

When the explicit expressions of $f_2$ and $f_3$ given by \eqref{eq:AWr} are used, these recurrence relations become the ones of the Askey--Wilson or $q$-Racah polynomials.
For  example, one finds the $q$-Racah polynomials (see Section 3.2 in \cite{KLS}) for the algebra $\uc$ and for $x=0,\dots,2j$, 
 \begin{equation}
  p_n=\frac{q^{2jn}}{(q-q^{-1})^n} \frac{(bq^{-2j+1},cq^{-2j+1},q^{-4j};q^2)_n}{(-aq^{2+2n-4j};q^2)_n}   R_n\Big(\mu(x),b q^{-2j-1}, -\frac{a}{b}q^{-2j-1}, q^{-2-4j}, -\frac{bc}{a} \Big| q^2\Big) \ ,
 \end{equation}
 with 
 \begin{equation}
  \mu(x)=q^{-2x}-\frac{bc}{a}q^{2x-4j}  , \qquad  \lambda=\frac{1}{q-q^{-1}}\left(q^{2j}\mu(x)-2\frac{bc}{a}q^{-2j} \right) .
 \end{equation}

\subsection{Specializations of the Askey--Wilson algebra}

We here turn to specializations or contractions of the realization of the Askey--Wilson algebra given by \eqref{eq:AWr}.
We provide different examples where the recurrence relations of the polynomials are recovered.
The number of free parameters decreases as the algebra is more and more specialized.

\paragraph{Specializations of the general results.} 

If we set the parameter $c=0$ in the realization \eqref{eq:AWr} of the Askey--Wilson algebra, we obtain a simpler algebra and the recurrence relations
are those of the $q$-Hahn polynomials or of the big $q$-Jacobi. If we set $b=c=0$, we get the recurrence relations for $q$-Krawtchouk or little $q$-Jacobi. 
For example, for $\uc$ and $a \in  {\rm \mathbb{R}}$, one gets the q-Krawtchouk polynomials (see Section 3.15 in \cite{KLS})
 \begin{equation}
  p_n =\frac{(-1)^n q^{2jn}(q^{-4j};q^{-2})_n}{(q-q^{-1})^n (-aq^{-4j+2n};q^{2})_n}     K_n(q^{-2x} ;aq^{-4j} ,2j|q^{2}),  \quad 
  \lambda=- \frac{q^{2(j-x)}}{q-q^{-1}}   \quad \text{for } x=0,1,\dots 2j.
 \end{equation}

\paragraph{Dual q-Hahn algebra.}

The operators with the upper (resp. lower) sign for $\uc$ (resp. $\unc$),
\begin{eqnarray}
\cX&=& K^{-2}\ , \\
\cY&=&E+ K^2 \frac{(q+q^{-1})K^2- {C}+\mu+\nu}{q-q^{-1}} 
 \mp q{K}^{2} \left( q {K}^{2}+\mu \right)  \left( q{K}^{2} +\nu \right) F\ ,
\end{eqnarray}
satisfy the following specializations of the $AW(3)$ relations:
\begin{eqnarray}
{} [\cX,[\cX,\cY]_q]_{q^{-1}} &=&  (q-q^{-1})(  C-\mu-\nu) \cX -q^2+q^{-2}, \label{eq:AW3ff} \\
{}[\cY,[\cY,\cX]_q]_{q^{-1}} &=& (q-q^{-1})(  C-\mu-\nu)  \cY -\mu\nu(q+q^{-1})^2 \cX +(q+q^{-1})(\mu\nu {C}-\mu-\nu).\label{eq:AW4ff} 
\end{eqnarray}

When $\mu\nu\neq 0$, the recurrence relations are those of the dual $q$-Hahn polynomials or the continuous dual $q$-Hahn polynomials.  
Explicitly, for $\uc$, one gets the dual $q$-Hahn polynomials (see Section 3.7 in \cite{KLS})
 \begin{eqnarray}
 && p_n =\frac{\mu^n \ (-q^{2j-1}/\mu,q^{4j};q^{-2})_n}{q^{2jn}(q-q^{-1})^n} \     R_n\left( q^{2x}+\frac{\nu}{\mu}q^{4j-2x};-\frac{ q^{2j+1}}{\mu},-\nu q^{2j+1},2j \Big|q^{-2}\right),  \\
  && \lambda= \frac{ \mu q^{2(x-j)}+\nu q^{2(j-x)} }{q-q^{-1}}   \quad \text{for } x=0,1,\dots 2j \ .
 \end{eqnarray}

If $\nu= 0$ (we can also choose $\mu=0$), one finds the recurrence relations of the affine $q$-Krawtchouk, quantum $q$-Krawtchouk, big $q$-Laguerre and $q$-Meixner polynomials.
For example, for $\uc$ and
 \begin{itemize}
  \item  for $\mu<-q^{-2j+1}$, one obtains the quantum $q$-Krawtchouk polynomials (see Section 3.14 in \cite{KLS}) and
 \begin{equation}
  p_n =\frac{q^{4jn+n}(q^{-4j};q^{2})_n}{(q-q^{-1})^n q^{2n^2}}     K^{qtm}_n(q^{-2x};-\mu q^{-2j-1},2j|q^{2}),  
  \quad \lambda= \frac{\mu q^{2(j-x)} }{q-q^{-1}}   \quad \text{for } x=0,1,\dots 2j\ ;
 \end{equation}
  \item  for $\mu<-q^{2j-1}$, one has the affine $q$-Krawtchouk polynomials (see Section 3.16 in \cite{KLS}) and
 \begin{equation}
  p_n =\frac{\mu^n q^{2jn}(-q^{2j-1}/\mu,q^{4j};q^{-2})_n}{(q-q^{-1})^n }     K^{aff}_n(q^{2x};-\frac{q^{2j+1}}{\mu} ,2j|q^{-2}),  
  \quad \lambda= \frac{\mu q^{2(x-j)} }{q-q^{-1}}   \quad \text{for } x=0,1,\dots 2j\ .
 \end{equation}
 \end{itemize}

\paragraph{$q$-Lie type algebras.}

We set
\begin{eqnarray}
\cX&=& K^{-2} \ ,\\
\cY&=&E- \frac{ a}{q-q^{-1}} K^2 + q {K}^{2}F \ .
\end{eqnarray}
These $\cX$ and $\cY$ satisfy the following specializations of the $AW(3)$ relations
\begin{eqnarray}
{} [\cX,[\cX,\cY]_q]_{q^{-1}} &=& (q-q^{-1})a \cX , \label{eq:AWql3} \\
{}[\cY,[\cY,\cX]_q]_{q^{-1}} &=& (q-q^{-1}) \cY + (q+q^{-1})^2 \cX -(q+q^{-1})  C.\label{eq:AW4ql} 
\end{eqnarray}
Again, the eigenvectors diagonalizing $\cY$ are associated to various polynomials.
For example,
\begin{itemize}
 \item for $\uc$ and $ a \in  {\rm \mathbb{R}}$, one gets the dual q-Krawtchouk polynomials (see Section 3.17 in \cite{KLS}) 
 by letting $ a=c-1/c$ with $c \in \mathbb{R}$ and, for $x=0,1,\dots, 2j$,
 \begin{equation}
  p_n =\frac{(q^{4j};q^{-2})_n}{(q-q^{-1})^nq^{2jn}c^n}     K_n(q^{2x}-c^2q^{-2x+4j};-c^2,2j|q^{-2}),  \quad 
  \lambda= \frac{q^{2(x-j)}/c-cq^{2(j-x)}}{q-q^{-1}}  ;
 \end{equation}
\item for $\unc$, one obtains:
\begin{itemize}
\item  the Al-Salam Chihara polynomials (see Section 3.8 in \cite{KLS}), by defining $c$ by $q^{2\ell} c +\frac{1}{cq^{2\ell}}=-\epsilon a$ and $\epsilon  \in \{1,-1\} $,
 \begin{equation}
 p_n =\frac{\epsilon^n}{(q-q^{-1})^n} Q_n\left(x;c  ,\frac{1}{cq^{4\ell}} \Big| q^{-2}\right) ,   \quad 
 \lambda= \frac{ 2\epsilon x}{q-q^{-1}} \quad \text{for }x \in [-1,1].
 \end{equation}
 Let us remark that for $c>1$ (which is the condition such the orthogonality of the polynomials to hold), the parameter of the algebra satisfies $a \in  ]-\infty,-2]  \cup [2,\infty[$;
\item the $q$-Meixner Pollaczeck polynomials (see Section 3.9 in \cite{KLS}), where
$ a=-2 \epsilon \cos {\phi}$  with $\epsilon  \in \{1,-1\} $,
 \begin{equation}
  p_n = \frac{\epsilon^n (q^{-2};q^{-2})_n}{(q-q^{-1})^n}P_n(x;q^{-2\ell}|q^{-2}), \quad
  \lambda= \frac{2\epsilon x}{q-q^{-1}} \quad \text{for }x=\cos(\theta+\phi) \text{ and } \theta \in [-\pi,\pi].
  \end{equation}
  In this case, the parameter of the algebra satisfies $ a \in [-2,2]$.
\end{itemize}
\end{itemize}

\section{Conclusions \label{sec:conc}}

Summing up, we have provided a new realization of the Askey--Wilson algebra in terms of the quantum algebras $\uc$ and $\unc$.
A new realization for the Racah algebra in terms of Lie algebras was also given.
We have shown that the representations we have introduced allow to obtain directly the recurrence relations of the orthogonal polynomials at the top of the Askey scheme : the $(q-)$Racah or the Wilson polynomials.  
We have also specialized the general case to get different simpler algebras and showed that the recurrence relations of other polynomials of the Askey scheme can also be recovered.
The merit of this work is to show that it is possible to find an algebraic realizations of the Askey--Wilson algebra in terms of simpler algebras in a way that 
the recurrence relation of the associated polynomials are recovered directly. It would be interesting to use these new realizations to construct in particular some realizations of the algebraic Heun 
operators recently introduced \cite{GrVZ}. The case $q=-1$ would deserve also a careful investigation: it corresponds to the Bannai--Ito algebra and the eponym polynomials.
\\

\vspace{1cm}

\noindent  {\bf Acknowledgements:}
N.~Cramp\'e has gratefully held a CRM-Simons Professorship during the course
of this project.
The research of L.~Vinet is supported
in part by a Discovery Grant from the Natural Science and Engineering
Research Council (NSERC) of Canada.

\end{document}